\documentclass{article}
\usepackage{amssymb}
\newcommand{\eq}{\begin{equation}}
\newcommand{\en}{\end{equation}}
\newcommand{\te}{\rightarrow}
\newcommand{\Phit}{\widehat{\Phi}}

\newcommand{\LL}{{\rm L}}

\newcommand{\re}[1]{\mbox{(\ref{#1})}}
\newcommand{\rem}[1]{\mbox{\rm (\ref{#1})}}

\newcommand{\tOm}{\widetilde{\Om}} 
 
\newcommand{\Om}{\nu}

\def\endpf{\hfill $\Box$ \vskip0.5cm}

\newtheorem{theorem}{\large Theorem}

\newtheorem{lemma}[theorem]{\large Lemma}

\begin{document}

\title{Asymptotic laws for regenerative compositions:\\
gamma subordinators and the like 
\thanks{Research supported in part by N.S.F. Grant DMS-0071448}
}
\author{Alexander Gnedin\,\,\thanks{Utrecht University; e-mail gnedin@math.uu.nl}\,\,,
\hspace{.2cm}
Jim Pitman\,\,\thanks{University of California, Berkeley; e-mail pitman@stat.Berkeley.EDU} 
\hspace{.2cm}
and 
\hspace{.2cm}
Marc Yor\,\,\thanks{University of Paris VI}
\\
\\
}
\date{\today}
\maketitle

\centerline{\bf Abstract}
For $\widetilde{\cal R} = 1 - \exp(- {\cal R} )$ a 
random closed set obtained by exponential transformation of the closed
range ${\cal R}$ of a subordinator,
a regenerative composition of generic positive integer $n$ is defined by recording the sizes of clusters 
of $n$ uniform random points as they are
separated by the points of $\widetilde{\cal R}$.
We focus on the number of parts $K_n$ of the composition 
when $\widetilde{\cal R}$ is derived from a gamma subordinator. We prove logarithmic asymptotics 
of the moments and central limit theorems for $K_n$ and other functionals of the composition
such as the number of singletons, doubletons, etc.
This study complements our previous work on asymptotics of these functionals
when the tail of the L\'evy measure is regularly varying at $0+$.

\vskip0.5cm

{\it AMS} 2000 {\it subject
classifications.} Primary 60G09, 60C05.

\newpage


\section{Introduction}
\label{Introduction}

For a drift-free subordinator
$(S_t,t\geq 0)$ define $\widetilde{\cal R}\subset [0,1]$ to be the range 
of the {\it multiplicative subordinator} $\widetilde{S}_t=1-\exp(-S_t)\,,t\geq 0$, obtained by the exponential 
transform. The {\it gaps} in $\widetilde{\cal R}$ are the interval components of the 
open set $[0,1]\setminus \widetilde{\cal R}$.
The following construction of random compositions of integers has been studied in
\cite{ gnedin03bs, gnedin03tsf, gnedinp03, gpy03regvar, csp}.
For each $n$, let 
$u_1,\ldots, u_n$ be an independent sample from the uniform distribution on $[0,1]$,
also independent of $\widetilde{\cal R}$. Define 
an ordered partition of the sample into nonempty blocks by assigning
two sample points $u_i$ and $u_j$ to the same block if and only if the points fall in the same gap
in $\widetilde{\cal R}$.
A composition ${\cal C}_n$ of integer $n$ is defined by the sequence of sizes of blocks, ordered from 
left to right.
The sequence of compositions $({\cal C}_n)$ is a regenerative composition structure in the sense of
\cite{gnedinp03}. The regeneration property means that, conditionally given the first part of ${\cal C}_n$ is $r$,
the remaining composition of $n-r$ has the same distribution as ${\cal C}_{n-r}$.

\par The purpose of this paper is to prove central limit theorems and determine the asymptotic behaviour of moments
for the number of parts $K_n$ and some other functionals of the composition ${\cal C}_n$. We focus on the
gamma subordinator with parameter $\theta > 0$, that is with L{\'e}vy measure $\nu({\rm d}y)=(e^{-\theta y}/y)\,{\rm d}y\,,$.
More generally, we obtain similar results
for L{\'e}vy measures which are sufficiently like the
gamma L{\'e}vy measure in their asymptotic behaviour at both $0+$ and
$\infty -$.
In the case of gamma subordinators our principal result (Theorem \ref{clt})
specialises as follows:

\begin{theorem}\label{gammaK} For the gamma subordinator with parameter $\theta>0$ the number of parts $K_n$ is asymptotically normal,
with ${\mathbb E\,}K_n\sim (\theta/2)\log^2 n$ and ${\rm var\,} K_n\sim (\theta/3)\log^3 n$.
\end{theorem}
\vskip0.2cm
\noindent
Other variables under consideration are $K_n(t)$, the number of parts of a partial composition 
produced by the subordinator restricted to $[0,t]$, and the
counts
$K_{n,r}$ and $K_{n,r}(t)$ defined as the multiplicity of part $r$ in ${\cal C}_n$ and the multiplicity 
of part $r$ in the partial composition, respectively. 

\par In our previous paper \cite{gpy03regvar} we studied $K_n$ for
a subordinator whose L{\'e}vy measure is regularly varying
at $0$ with index $\alpha\in \,]0,1]$. In that case we obtained
limit theorems for $K_n$ with normalisation by
$n^\alpha\,\ell(n) \,$ where $\ell$ is a function of slow variation at 
$\infty$. 
According to these results, $K_n/(n^\alpha\,\ell(n))$
converges strongly and with all moments to a nondegenerate limit, which is not gaussian. Moreover, $K_n(t)$ 
is of the same order of magnitude as $K_n$, for each $t>0$. 
For the case of a compound Poisson process, when the L{\'e}vy measure is finite, 
Gnedin \cite{gnedin03bs} had previously shown 
that $K_n$ typically exhibits 
a normal limit with
both expectation and variance of the order $\log n$, though $K_n(t)$ remains bounded as $n$ grows, for each $t>0$.
Thus the case of gamma-type subordinators may be seen as intermediate between the above two:
as in the regular variation case, $K_n(t)$ is unbounded, but, as in the case of a
finite L{\'e}vy measure,
its contribution to $K_n$
is asymptotically negligible, for each $t>0$. 

\par Results akin to Theorem \ref{gammaK} are very different from the limit theorems in \cite{gnedin03bs, gpy03regvar}
and require other techniques. Our approach here
relies on a recent version of the contraction method
due to Neininger and R{\"u}schendorf \cite{neinruesch}. Application of this method requires an appropriate decomposition of $K_n$, 
control of the principal asymptotics of two moments, and an 
estimate of the remainder term in the asymptotic expansion of the variance.
We provide this background by poissonising ${\cal C}_n$ and
applying the
Mellin transform technique to analyse integral recursions for the moments. 
\par We also establish joint convergence to a gaussian limit for the sequence of small-part counts
$(K_{n,r}, r\geq 1)$, as $n\to\infty$.
This kind of convergence holds neither 
in the regular variation case (when the scaled $K_{n,r}$'s converge to multiples of 
the same variable) nor in the compound Poisson case 
(when $(K_{n,r}, n\geq 1)$ is bounded uniformly in $n$, for each fixed $r$). 
It resembles, however, Karlin's \cite{karlin67urn} central limit theorem 
for nonrandom frequencies in the case of regular variation with a positive index.

\section{Setup and notation}

Each {\it part} of ${\cal C}_n$ corresponds to a nonempty cluster of 
the uniform sample points within some gap. So $K_n$ is the number of gaps
occupied by at least one of $n$ sample points, and $K_{n,r}$ 
is the number of gaps that contain exactly $r$ sample points. 
Keep in mind that $K_n=\Sigma_r\,K_{n,r}, $ and $n=\Sigma_r\, r K_{n,r}$. Sometimes $K_n$ is called the {\it length},
and $n$ the {\it weight} of the composition.
Similarly,
$K_n(t)$ and $K_{n,r}(t)$ are the counts of sizes of clusters within the gaps of a 
smaller set $\widetilde{\cal R}\cap [0,\widetilde{S}_t]$, corresponding to the 
subordinator restricted to $[0,t]$.

\par 
For $\nu$ the L{\'e}vy measure of $(S_t)$ let $\tOm$ be the measure on $[0,1]$ obtained from $\nu$ by the 
exponential transform $y\mapsto 1-e^{-y}$. The 
Laplace exponent of the subordinator is given by the formulas
\begin{equation}\label{Lapex}
\Phi(s)=\int_0^\infty (1-e^{-sy})\nu({\rm d}y)=\int_0^{1}(1-( 1-x)^s)\tOm({\rm d}x)\,,\qquad \Re\, s\geq 0\,.
\end{equation}
Considering also the binomial moments 
\begin{equation}\label{2Lapex}
\Phi(n:m):=
\int_0^{1}{n\choose m}x^m (1-x)^{n-m}\tOm({\rm d}x)\,,\qquad 1\leq m\leq n\,,
\end{equation}
the ratio $\Phi(n:m)/\Phi(n)\,, m=1,\ldots,n,$ gives the distribution of the first part of ${\cal C}_n$,
and the probability of a particular value of ${\cal C}_n$ is a product of factors 
of this type \cite{gnedinp03}.
We also record the relation between power and logarithmic moments of measures $\nu$ and $\tOm$:
\begin{equation}\label{logmom}
{\tt m}_j:=\int_0^{\infty} y^j \,\nu({\rm d}y)=\int_0^1 |\log (1-x)|^j\,\tOm ({\rm d}x)\,,~j=1,2,\ldots\,.
\end{equation}



\par As in \cite{gpy03regvar}, we shall also consider the {\it poissonised}
composition $\widehat{\cal C}_\rho$ derived from $\widetilde{\cal R}$ 
by the same construction as ${\cal C}_n$, but 
with the set of atoms of
a homogeneous Poisson point process
on $[0,1]$ with rate $\rho$,
instead of the uniform sample of fixed size.
We denote $\widehat{K}_\rho,\, \widehat{K}_{\rho,r}, \,
\widehat{K}_\rho(t),\, \widehat{K}_{\rho,r}(t)$ the obvious counterparts of the fixed-$n$ variables, and
introduce the Poisson transform of the Laplace exponent
\begin{equation}\label{P-transform}
\Phit(s) :=\int_0^{1} (1-e^{-s x})\tOm({\rm d}x)\,
= \, e^{-s}\sum_{n=0}^{\infty} {s^n \over n!}\,\Phi(n) .
\end{equation}
We shall be using throughout the abbreviations

$${\rm L}=\log \rho,\qquad \,\log^a x=(\log x) ^a\,.$$ 

\vskip0.5cm

\par Our primary concern are the compositions induced by a {\it gamma subordinator}, whose 
L{\'e}vy measure, its exponential transform and the 
 Laplace exponent are given in terms of a parameter $\theta>0$ by
\begin{equation}\label{gamma}
\nu({\rm d}y)={e^{-\theta y}\over y}\,{\rm d}y\,,~y\in \,]0,\infty[\,,\quad
\tOm({\rm d}x)={(1-x)^{\theta-1}\over |\log (1-x)|}\,{\rm d}x\,,~x\in \,]0,1]\,.~~
\quad \Phi(s)=\log\left(1+{s\over \theta}\right)\,,
\end{equation}
The density 
of $\tOm$ has a pole and the tail $\tOm[x,1]$ has a logarithmic singularity at $0$. 
We could consider a larger family of measures which differ from (\ref{gamma}) by a positive factor,
but this would not increase generality, since multiplying the L{\'e}vy measure by such a factor does not affect
the laws of ${\cal R}$, ${\cal C}_n$, $K_n$ or $\widehat{K}_{\rho}$ (although it requires a linear 
time-change for functionals of
a partial composition like
$K_n(t)$ or $\widehat{K}_{\rho}(t)$). 
\vskip0.5cm

\par More generally, we consider a drift-free subordinator with 
L{\'e}vy measure $\nu$, which has a continuous density on $]0,\infty[$ and
satisfies  the following conditions
(L) and (R)
which we record in terms of $\nu$ and $\tOm$: 

\begin{enumerate}
\item[(L)] either of the following four equivalent conditions holds
$$\Phi(\rho)= {\rm L}+c+O(\rho^{-\epsilon}) \,,{\rm ~as~} \rho\to\infty\,,\qquad \Phit(\rho)= {\rm L}+c+
O(\rho^{-\epsilon})\,, {\rm ~as~} \rho\to\infty$$
$$\nu[y,\infty]= - \log y +c-\gamma+ O(y^{\epsilon})\,, {\rm ~as~} y\downarrow 0\,,
\qquad \tOm[x,1]= - \log x+c-\gamma +O(x^{\epsilon})\,, {\rm ~as~} x\downarrow 0\,$$ 
where $c$ and $\epsilon>0$ are some constants and 
$\gamma=-\Gamma'(1)$ is the Euler constant, 
\item[(R)] either of the following two equivalent conditions holds
$$\nu[y,\infty]=O(e^{-\epsilon y})\,, ~y\uparrow \infty\,,\qquad \tOm[x,1]=O((1-x)^{\epsilon})\,,~ x\uparrow 1$$
where $\epsilon>0$.
\end{enumerate}
See Appendix for the equivalence of conditions in (L) 
(L stands for {\it left} and {\it logarithmic} at $0$).
Condition (R) (R for {\it right} or {\it regular at $1$}) implies that all moments ${\tt m}_j$ are finite,
and that $\Phi$ is analytical for $\Re\, s>-\epsilon$.

\par Throughout in this paper, $\epsilon$ or $\delta$ are some sufficiently small positive constants whose values
 are context-related
and may change from line to line. We denote $c$ the constant in (L).
For a generic positive constant  we write  $d$,
and use $c_j,\, d_j$  for further real constants which are not important, but may appear with explicit evaluation
in asymptotic expansions or other formulas.
For shorthand, the right tail of $\tOm$ is denoted by
$$\vec{\nu}(x)=\tOm[x,1]\,,\qquad x\in \,]0,1]$$
and a homogeneous Poisson point process with intensity $\rho$ on $[0,1]$ is denoted PPP($\rho$).

\section{Recursions}
We shall make use of two types of recursions. These apply
to regenerative compositions generated by a drift-free
subordinator, without any special assumptions on the L{\'e}vy measure.
(They can also be readily generalised to the case with drift). 
The first type of recursion, suited to the Poisson framework, is based
on a decomposition of $\widetilde{\cal R}$
analogous to the first-jump decomposition of a renewal process.
The second type is based on splitting the range of a subordinator 
by its value at a stopping time.
\vskip0.5cm

\noindent
{\bf Integral recursions} 
Define a {\it pattern} $E$ to be a nonempty set of positive integers.
We say that a finite configuration of points within a given subinterval $]a,b[\,\in \,]0,1[\,$ 
fits in $E$ if the cardinality
of the configuration is in $E$.
For $E$ fixed, let 
$$\pi(\rho)=\sum_{r\in E} e^{-\rho} \rho^r/r!$$ 
be the probability that the configuration of atoms of 
PPP($\rho$) fits in the pattern $E$.
The probability that the PPP($\rho$) configuration on a subinterval $]a,b[$ fits in $E$ is equal to $\pi(\rho(b-a))$.

\par Consider the 
poissonised composition induced by an arbitrary 
drift-free subordinator with $\tOm\{1\}=0$ (no killing). 
For $E$ fixed, let $N_{\rho}$ be the number of gaps of $\widetilde{\cal R}$ such that the PPP($\rho$)
configuration within the gaps fits in $E$. 
The count $N_{\rho}$ is a functional of the poissonised 
composition $\widehat{\cal C}_\rho$; in particular, defining $E$ to be $\{1,2,\ldots\}$
(the configuration is nonempty) or $E=\{r\}$ (the configuration consists of $r$
points) we obtain $N_{\rho}=\widehat{K}_{\rho}$, respectively $N_{\rho}=\widehat{K}_{\rho,r}$.
Other possibilities may be considered, for example taking 
$E=\{1,3,5,\ldots\}$
the count $N_{\rho}$ becomes the number of odd parts of the composition.

\par Let $p_j(\rho)={\mathbb P}(N_{\rho}=j)$ be the distribution of $N_{\rho}$ for some fixed pattern.
Each $p_j$ may be extended to an entire function
of a complex variable, with the initial value 
$p_j(0)=1(j=0)$. 
Introduce the factorial moments
$$f^{(m)}(\rho)={\mathbb E}N_{\rho}(N_{\rho}-1)\cdots (N_{\rho}-m+1)\,,\,\,\,m=0,1,\ldots$$
with the convention
$f^{(0)}(\rho)=1$. The following lemma is a minor variation of \cite[Lemma 6.1]{gpy03regvar}.
\vskip0.5cm

\begin{lemma}
\label{pattern}
Let $E$ be a pattern with the probability of occurrence $\pi(\rho)$. 
The distribution of $N_{\rho}$
satisfies the integral recursion
\begin{equation}\label{p-rec}
\int_0^1 (p_j(\rho) -(1-\pi(\rho x))\, p_j(\rho(1-x))\,\tOm({\rm d}x)\,=\int_0^1 \pi(\rho x)p_{j-1}(\rho(1-x))\, 
\tOm({\rm d}x)
\end{equation}
for $j=1,2,\ldots$ and $\rho\geq 0$. The same equation is also valid for $j=0$, with $0$ on the right-hand side.
Furthermore, the factorial moments of $N_\rho$ satisfy the recursion
\begin{equation}\label{f-rec}
\int_0^1 (f^{(m)}(\rho) -f^{(m)}(\rho(1-x))\,\tOm({\rm d}x)\,=m\, \int_0^1 \pi(\rho x)f^{(m-1)}(\rho(1-x))\,
\tOm({\rm d}x)
\end{equation}
which taken together with $f^{(0)}(\rho)=1$ and $f^{(m)}(0)=1(m=0)$ 
uniquely determines them. 
\end{lemma}
\noindent
{\it Proof.} The derivation of recursions follows as in \cite[Lemma 6.1]{gpy03regvar}. 
The uniqueness claim is a consequence of analyticity.
\endpf

\vskip0.5cm
\par Integral recursions for the distribution and the factorial 
moments of $\widehat{K}_{\rho}$ follow by taking $\pi(\rho)=1-e^{-\rho}$ which is the probability that 
the PPP($\rho$) configuration is nonempty; we have then
$p_0(\rho)=e^{-\rho}$.
To 
obtain recursions for $\widehat{K}_{r,\rho}$ we should take $\pi(\rho)=e^{-\rho}\rho^r/r!$ (in this case no simple formula for 
$p_0(\rho)$ is known).
\vskip0.5cm
\noindent
{\bf Splitting at an independent exponential time}
Further recursions follow by
splitting the range of a subordinator by its value at a stopping time. 
Though the fixed-$n$ version is needed to apply \cite{neinruesch}, we focus on the poissonised model,
which simplifies moment computations.
Transfer of results to the fixed $n$ model then follows by elementary depoissonisation.

\par For each $t\geq 0$, consider the sigma-algebra generated by both
$(\widetilde{S}_u, u\leq t)$ and the restriction of PPP($\rho$) to $[0,\widetilde{S}_t]$.
As $t$ varies, this defines a filtration, and we can consider
stopping times with respect to it.
Let $\tau$ be such a stopping time 
with range $[0,\infty]$
and let $b_\rho=\widehat{K}_\rho(\tau)$ be the number of parts of ${\cal C}_\rho$ produced by the subordinator up to time 
$\tau$.
The strong Markov property of $\widetilde{\cal R}$ 
along with the independence property of PPP
imply that the
number of parts 
satisfies the distributional equation
\begin{equation}\label{fpoint}
\widehat{K}_\rho \stackrel{d}{=} b_\rho+\widehat{K}_{I_\rho}'
\end{equation}
where $b_\rho=\widehat{K}_\rho(\tau)\,, I_\rho=\rho (1-\widetilde{S}_\tau)$ 
and $(\widehat{K}_\rho')$ is a distributional copy of $(\widehat{K}_\rho)$, independent of $(b_\rho, I_\rho)$. 
For example, letting $\tau$ be the first time that a jump of the subordinator covers at least one Poisson point,
the equation holds with $b_\rho=1$.
Identities analogous to (\ref{fpoint}) can be written also for $\widehat{K}_{\rho,r}$ and more general pattern counts.
\vskip0.5cm

\par We shall consider in some detail the most important 
case when $\tau$ is an exponential time with rate $\lambda$, independent of the subordinator and the PPP
(thus $\tau$ is a {\it randomised} stopping time with respect to the above filtration).
The stopped process has an obvious interpretation as a {\it killed} multiplicative
subordinator, which jumps at time $\tau$ to the terminal value $1$ (thus producing the final {\it meander} gap 
in the range).
In accord with our previous notation, let $\widehat{K}_{\rho}(\tau)$ be the number of parts produced by the subordinator {\it before} $\tau$ 
(thus a possible
part induced by
the meander gap is not taken into account). Let $(p_j(\rho), j\geq 0)$ and $(f^{(m)}, m\geq 0)$ be the distribution
and factorial moments of $\widehat{K}_{\rho}(\tau)$.

\begin{lemma}
\label{split}
For $\tau$ an exponential time with rate $\lambda$,
independent of the process $(\widehat{K}_{\rho}(t), t \ge 0)$,
the factorial moments of $\widehat{K}_{\rho}(\tau)$ satisfy the recursion
\begin{equation}\label{f-split-rec}
\lambda \,f^{(m)}(\rho)+\int_0^1 (f^{(m)}(\rho) -f^{(m)}(\rho(1-x)))\,\tOm({\rm d}x)\,=m\, \int_0^1 \pi(\rho x)f^{(m-1)}(\rho(1-x))\,
\tOm({\rm d}x)
\end{equation}
which taken together with $f^{(0)}(\rho)=1$ and $f^{(m)}(\rho)=1(m=0)$
uniquely determines them.
\end{lemma}
{\it Proof.} In the case of a finite L{\'e}vy measure,
the first-jump decomposition of the range yields
\begin{eqnarray*}
p_j(\rho)=
\int_0^1 \bigg((1-\pi(\rho x)) p_j(\rho(1-x))+\pi(\rho x)\, p_{j-1}(\rho(1-x))\bigg)\,{\tOm({\rm d}x)\over
\vec{\nu}(0)+\lambda}\\
p_0(\rho)= {\lambda\over \vec{\nu}(0)+\lambda}+
\int_0^1 (1-\pi(\rho x))p_{0}(\rho(1-x))
{\tOm({\rm d}x)\over
\vec{\nu}(0)+\lambda}\,.
\end{eqnarray*}
To see the extension for an arbitrary L{\'e}vy measure, substitute
$$p_j(\rho)=p_j(\rho){\lambda\over \vec{\nu}(0)+\lambda}+
\int_0^1 p_j(\rho) 
{\tOm({\rm d}x)\over
\vec{\nu}(0)+\lambda}\,,$$
then rearrange terms and argue as in \cite[Lemma 6.1]{gpy03regvar}.
\endpf

\vskip0.5cm
\noindent
If the meander gap is counted (this corresponds to $\widehat{K}_\rho$ for killed subordinator) we should modify the
recursion for moments by adding $\lambda \pi(\rho)1(m=1)$ to the right-hand side.
\vskip0.5cm
\noindent
{\bf Remark.} The distribution and moments of $\widehat{K}_\rho(t)$ may be obtained as the inverse Laplace transform 
in $\lambda$
from the analogous quantities for $\widehat{K}_\rho(\tau)$.

\section{Moments of $\widehat{K}_\rho$}
\label{momKrho}
\noindent
{\bf Mellin transform resolution}
The recursion (\ref{f-rec}) is intrinsically related to a convolution-type integral equation
\begin{equation}\label{int-eq}
\int_0^1 (f(\rho)-f(\rho(1-x)))\tOm({\rm d}x)=g(\rho)
\end{equation}
with the function $g$ given and the function $f$ unknown. For $g$ analytic, there is a unique analytic solution $f$ with 
given initial value $f(0)$. Observe that constants $d$ constitute the null space of the integral operator, thus
for each solution $f$ the function $f+d$ is another solution, with the initial value $f(0)+d$.

\par It is easy to write out a power series solution to (\ref{int-eq}), but such a representation itself does not help
describe the large-$\rho$ behaviour. We turn therefore to asymptotic methods based on the Mellin transform.
Recall that for a locally integrable real-valued function $\phi$ on $[0,\infty]$ the Mellin transform is defined by the integral
$$M\phi(s)=\int_0^{\infty} \rho^{s-1}\phi(\rho) \,{\rm d}\rho\,$$
which is assumed to converge absolutely for $s$ in some open interval on the real axis, 
hence also for $s$ in the open vertical strip based on the interval. 
The left convergence abscissa of $M\phi$ is determined by the behaviour of $\phi$ near $0$, while the right 
convergence abscissa is determined by the behaviour of $\phi$ at $\infty$.

\par The  analysis to follow is based on the formula 

\begin{equation}\label{Mf}
M f(s)={M g(s)\over \Phi(-s)}
\end{equation}
which is valid in the common domain of definition of all ingredients.
The formula follows by Fubini and a change of variable from 
$$\int_0^\infty \rho^{s-1}{\rm d}\rho \int_0^1(f(\rho)-f(\rho(1-x))\tOm({\rm d}x)=Mf(s)
\int_0^1(1-(1-x)^{-s})\tOm({\rm d}x)=Mf(s)\Phi(-s)\,$$
according to formula \re{Lapex}.
Another important tool is the following correspondence between 
asymptotic expansion of a function and singularities of its Mellin transform. 
See \cite{flaj-harm} for a fuller exposition of this technique.

\begin{lemma} 
\label{singan}
{\em \cite[Section 2]{flaj-harm} }
Suppose the Mellin transform 
 $Mf$ of a function $f$ is analytic in a strip $a<\Re s <b$.
If $Mf$
can
be extended meromorphically 
through the right convergence abscissa 
in a larger trip $a<\Re s<b+\epsilon$ 
 and has finitely many poles there,
then 
each pole $z$ and 
each term $d\,(s-z)^{-k-1}$ in the Laurent expansion of $M f$ at $z$ contributes 
the term $${d \,(-1)^{k+1}\over k!}\,\,\rho^{-z}\,\log^k \rho$$ to the asymptotic expansion of 
$f$ at $\infty$. The remainder term of the expansion of $f$ is then $O(\rho^{-b-\epsilon})$.
Conversely, the asymptotic expansion of $f$ at $\infty$ with such a remainder 
implies the termwise singular expansion 
of $Mf$ in the strip $a<\Re s<b+\epsilon$ provided that for some $\delta>0$
$$|Mf(s)|=O(|s|^{-1-\delta})$$
as $|s|\to\infty$ in a strip  $b'<\Re\, s<b+\epsilon$ for some $b'<b$.
\end{lemma} 
\vskip0.5cm
\noindent
\par For
the expectation $f^{(1)}(\rho):={\mathbb E}\,\widehat{K}_{\rho}$ 
the formulas (\ref{f-rec}) and (\ref{Mf}) specialise as 
\eq
\label{mf1s}
M f^{(1)}(s)={M \Phit(s)\over \Phi(-s)}
\en
where $\Phit(s)$ is defined by \re{P-transform}.
Using Fubini and integration by parts we transform the numerator as
\begin{eqnarray*} M \Phit(s)=\int_0^{\infty} \rho^{s-1} {\rm d}\rho \int_0^1 \rho e^{-\rho x}
\vec{\nu}(x)\,{\rm d}x=
\int_0^1 \vec{\nu}(x)\,{\rm d}x
\int_0^{\infty} \rho^{s} e^{-\rho x} {\rm d}\,\rho=\\
\Gamma(s+1)\int_0^1 x^{-s-1}\vec{\nu}(x){\rm d}x = -\Gamma(s)\int_0^1 x^{-s}\,\tOm({\rm d}x) .
\end{eqnarray*}
We can now re-write formula \re{mf1s} as
\begin{equation}\label{Mf1}
M f^{(1)}(s)=-\Gamma(s){\Phi(-s:-s)\over\Phi(-s)}
=-\Gamma(s)
{\int_0^1 x^{-s-1}\,\vec{\nu}(x)\,{\rm d}x\over \int_0^1 (1-x)^{-s-1}\,\vec{\nu}(x)\,{\rm d}x}\,,
\end{equation}
where we used (\ref{Lapex}) and (\ref{2Lapex}) in the form
$$\Phi(-s:-s)= -s\,\int_0^1 x^{-s-1}\,\vec{\nu}(x) 
\,{\rm d}x
\,,\qquad \Phi(-s)=-s\int_0^1 (1-x)^{-s-1}\,\vec{\nu}(x)
\,{\rm d}x
.
$$
\vskip0.5cm
\noindent
{\bf Expectation} Assuming (L)  we conclude from (\ref{Mf1}) that $Mf^{(1)}$ is defined in the strip $-1<\Re\, s<0$.
We focus therefore on the meromorphic continuation of $Mf^{(1)}$ through the right convergence abscissa $\Re s=0$.
The same formula says that  the poles of $Mf^{(1)}$ might be caused either by poles of $\Gamma$ 
or by poles of $\Phi(-s:-s)/s$, 
or by zeros of $\Phi(-s)/s$. 
We shall consider the three ingredients separately.

\par For $\Re s>-1$ the gamma function has a unique pole at $0$, with Laurent expansion 
$$\Gamma(s) = {s^{-1}}-\gamma+d_1 s + O(s^2)$$
where $d_1={\gamma^2/2}+{\pi^2/12}$. 
With reference to Stirling's approximation, when $u$ is a bounded real number and $v$ is a large real number
\begin{equation}\label{stirling}
|\Gamma(u+{\tt i}v)|\sim (2\pi)^{1/2}|v|^{u-1/2} e^{-\pi |v|/2}
\end{equation}
uniformly in $u$.

\par By assumption (R), the function
$$-\Phi(-s)/s= \int_0^1 (1-x)^{-s-1}\,\vec{\nu}(x)\,{\rm d}x\,$$
is analytic for $\Re s<\epsilon$, and its behaviour at complex infinity is regulated by the next lemma.
\vskip0.5cm
\begin{lemma} \label{fourier} If $\nu$ has a continuous density on $]0,\infty[$ and satisfies {\rm (L)} and {\rm (R)} then
$$\Phi(s)\sim \gamma\, \log|s|\,,\quad {\rm as~}|s|\to\infty\,$$
uniformly in any strip $-\epsilon<\Re\,s <d$, for $\epsilon$  sufficiently small.
\end{lemma}
{\it Proof.} The measure  $y\,\nu({\rm d}y)$ is finite, with exponential decay at infinity, thus
for $s=u+{\tt i}v$ the partial derivative
$${{\rm d}\over {\rm d}u}\Phi(u+{\tt i}v)=\int_0^\infty e^{-uy}e^{-{\tt i}vy}y\,\nu({\rm d}y)$$
is bounded for $u>-\epsilon$, uniformly in $v$. Hence  $|\Phi(u+{\tt i}v)-\Phi({\tt i}v)|$ is uniformly bounded both  in $v$ and 
in  $u\in\,[-\epsilon, d]$. Thus it is sufficient to show the asymptotics for $u=0$.
But   for the Fourier integral 
$$\Phi({\rm i}v)={\rm i}v\,\int_0^\infty e^{-{\rm i}vy}\nu[y,\infty]\,{\rm d}y$$
the desired asymptotics follows from the expansion (L) and smoothness of $\nu[\,\cdot\,,\infty]$ by application of
\cite[Section 3.1, Theorem 1.11]{Fedoryuk} (with a further reference to 
\cite{Armstrong}).\endpf

\noindent
By one further elementary lemma (see Appendix), condition (L) implies that $\Phi(-s)/s$ has no zeroes for $\Re\,s\leq 0$,
hence by Lemma \ref{fourier}   there are no zeroes in the strip $-1<\Re\,s\leq \epsilon$.  
The Taylor expansion at $0$ 
involves the logarithmic moments (\ref{logmom}) 
\eq
\label{phiss}
{-\Phi(-s)\over s} = {\tt m}_1 +{\tt m}_2 {s \over 2!}+{\tt m}_3{s^2\over 3!} + O(s^3). 
\en

\par The integral 
$${-\Phi(-s:-s)/ s}=\int_0^1 x^{-s-1}\vec{\nu}(x)\,{\rm d}x$$
converges for ${\Re }\,s<0$ and by 
assumption (L) 
this function is meromorphic for ${\Re}\,s<\epsilon$ and
in this half-plane has the unique pole at $s=0$, where the Laurent series starts with
$${-\Phi(-s:-s)\over s} = {s^{-2}}-{(c-\gamma)s^{-1}}+d_2 + o(1)\,$$
where $d_2$ is some constant.
The  meromorphic extension  is obtained from this formula, provided we
substitute a suitable analytic function for the $o(1)$ term.
This follows by computing 
$$\int_0^1 (|\log x|+c-\gamma+h(x))x^{-s-1}{\rm d}x\,$$
with $h(x)=O(x^\epsilon)$ and $\Re\,s<0$, and noting that
the integral of $h$ is analytic  for $\Re\,s<\epsilon$ and it is bounded in this domain due to a maximum 
ridge on the real axis.

\par It follows that $Mf^{(1)}$ is meromorphic in the strip $-1<\Re s<\epsilon$ 
with a unique singularity at $0$, where it has
a {\it triple} pole.
Putting the ingredients together and with a minor assistance of {\tt Mathematica} we obtain
the Laurent expansion 

$$M f^{(1)}(s) = -{1\over {\tt m}_1} s^{-3}+
\left({c\over{\tt m}_1}+{{\tt m}_2\over 2{\tt m}_1^2}\right){s^{-2}}+
d_3 {s^{-1}} + O(1)$$ 
where 
$$d_3=
\left({-d_1-d_2-c\gamma+\gamma^2\over{\tt m}_1}-{c{\tt m}_2\over 2{\tt m}_1^2}-
{{\tt m}_2^2\over 4{\tt m}_1^3}+{{\tt m}_3\over 6{\tt m_1}^2}\right)$$
is a constant whose explicit value will not be used below.
By Lemma \ref{singan}, 
as $\rho \te \infty$,
the asymptotic expansion
of the expectation is
\begin{equation}\label{asympexp}
f^{(1)}(\rho) = {1\over 2{\tt m}_1} {\rm L}^2 +\left({{\tt m}_2\over 2{\tt m}_1^2} +{c\over{\tt m}_1}\right){\rm L}
-d_3 + O(\rho^{-\epsilon}).
\end{equation}
The decay condition on  $Mf^{(1)}$  required in Lemma \ref{singan}
 is satisfied, because this $Mf^{(1)}$ has exponential decay as $|s|\to\infty$ in a strip about the imaginary axis
due to (\ref{stirling}), Lemma \ref{fourier} 
and because $|\Phi(s:s)|=O(|s|^{-1})$.

\vskip0.5cm

\noindent
{\bf Evaluating the right-hand side of (\ref{f-rec}).}
We will use (\ref{Mf}) once again, to compute $f^{(2)}$. 
The right-hand side of (\ref{f-rec})
is evaluated with the help of the next lemma.

\vskip0.5cm

\begin{lemma}\label{le2} Assume {\rm (L)} and {\rm (R)} and suppose $f$ is an increasing positive function
of $\rho \in [0,\infty]$, which admits an
asymptotic expansion at $\infty$ as a polynomial in $\,\rm L$, with 
a remainder $O(\rho^{-\epsilon})$. Then a similar expansion holds also for 
the function 
$$h(\rho) =\int_0^1 (1-e^{-\rho x})f(\rho(1-x))\tOm({\rm d}x)$$
and this expansion is obtained from that of $f$ by replacing each ${\rm L}^k$ with
$${\rm L}^{k+1}+c{\rm L}^k+\sum_{j=1}^{k} (-1)^j {k\choose j} {\tt m}_{j}{\rm\, L}^{k-j}\,$$
where $c$ is as in {\rm (L)}.
\end{lemma}
{\it Proof.} By the assumption
\begin{equation}\label{fpol}
f(\rho)=\sum_{j=0}^k d_j\,\LL^j+O(\rho^{-\epsilon}).
\end{equation}
Start with the case $f=\LL^k$. The integral is computed 
by the binomial expansion of 
$({\rm L}+\log(1-x))^k$ and termwise integration.
The leading term 
does not depend on $x$ and, in view of (L), it contributes ${\rm L}^{k}\,({\rm L}+c)+ O(\rho^{-\epsilon})$.
The lower order
terms are integrated using (\ref{logmom}), with a remainder estimated by a constant multiple of
$${\rm L}^k \,\int_0^1 e^{-\rho x} |\log (1-x)|^k\, \tOm({\rm d}x)\,,\quad k\geq 1.$$
The last integral is evaluated as $O(\rho^{1-\epsilon})$ using integration by parts, using (L) and applying
a standard  Tauberian theorem on Laplace transforms.

\par By linearity we obtain asymptotic expansion for the 
{\it model integral} $I(\rho)$ which corresponds to the polynomial part of $f$,
as in (\ref{fpol}) but with zero remainder term. 
One checks easily using (R) that for $\delta:=\rho^{-1/2}$
the contribution to $I(\rho)$ of the integral over $[1-\delta,\,1]$ is estimated as $O(\rho^{-\epsilon})$.

\par For a general $f$ as in (\ref{fpol}) we again split the integral at $1-\delta$.
For $x\in [0,1-\delta]$ we can apply expansion (\ref{fpol}) to $f(\rho(1-x))$ with a remainder $O(\rho^{-\epsilon})$, thus
by (L) we have
$$\int_0^{1-\delta}(1-e^{-\rho x})f(\rho(1-x))\tOm({\rm d}x)=I(\rho)+O(\rho^{-\epsilon}).$$
It remains to show that the integral over $[1-\delta,\,1]$ is $O(\rho^{-\epsilon})$ but this is easy because
the assumed monotonicity yields the bound 
$$\int_{1-\delta}^1 (1-e^{-\rho x})f(\rho(1-x))\tOm({\rm d}x)<f(\rho\,\delta)\,\vec{\nu}(1-\delta)$$
and this has the desired order by (\ref{fpol}), assumption (R) and our choice of $\delta$.
\endpf
\vskip0.5cm

\par Taking $f^{(1)}$ for $f$, the lemma translates the expansion (\ref{asympexp}) into the $\rho\to\infty$ 
formula
for the
right-hand side of (\ref{f-rec}) with $m=2$
\begin{equation}\label{2g}
g_1(\rho)= {1\over {\tt m}_1}{\rm L}^3+\left({3c\over {\tt m}_1}+{{\tt m}_2\over {\tt m}_1^2}\right){\rm L}^2
+d_4{\rm L}+d_5+O(\rho^{-\epsilon})
\end{equation}
where
$d_4, d_5$ are some constants,
and
the factor $2$ in (\ref{f-rec}) is included in the definition of $g_1$.

\vskip0.5cm
\noindent
{\bf Moments of the second order} 
Using the direct correspondence in Lemma \ref{singan}, we see that $0$ is the sole singularity
of $Mg_1$, and derive the expansion at $s=0$
$$Mg_1(s) = {6\over {\tt m_1}}s^{-4}-\left({6c\over {\tt m}_1}+{2{\tt m_2}\over {\tt m}_1^2}\right) s^{-3}+  
O(s^{-2})\,$$
which together with \re{Mf} and \re{phiss} implies
$$Mf^{(2)}(s)={Mg_1(s)\over \Phi(-s)} = - {6\over {\tt m}_1^2}\,s^{-5}+{6c{\tt m}_1+5{\tt m}_2\over {\tt m}_1^3}
\,s^{-4}
+ O(s^{-3})\,$$
Finally, the inverse correspondence applied to $f^{(2)}$ yields
\begin{equation}\label{f2as}
f^{(2)} (\rho) = {1\over 4{\tt m}_1^2}{\rm L}^4+\left( {c\over {\tt m}_1^2}+{5{\tt m}_2\over 6{\tt m}_1^3}\right){\rm L}^3 +
d\, {\rm L}^2 + O(\LL)
\end{equation}
provided the bound $|Mf^{(2)}(s)|=O(|s|^{-1-\epsilon})$ can be established.
This requires some effort because for $Mg_1$ 
no explicit formula is available. 
Postponing justification of the decay condition, and recalling 
$${\mathbb E}\,\widehat{K}_{\rho}=f^{(1)}(\rho)\,,~~{\rm var\,}\widehat{K}_{\rho}= f^{(2)}(\rho)+ f^{(1)}(\rho)- (f^{(1)}(\rho))^2$$
we have from (\ref{asympexp}) and (\ref{f2as})

\vskip0.5cm
\begin{theorem}
\label{momenty} If the L{\'e}vy measure of subordinator has a continuous density on $\,]0,\infty[$ 
and satisfies the assumptions {\rm (L)} and {\rm (R)}, then  
as $\rho \to \infty$ with ${\rm L } = \log \rho$,
the asymptotic expansions of the first two central moments
of $\widehat{K}_{\rho}$
are
\begin{equation}\label{f1as}
{\mathbb E}\,\widehat{K}_{\rho}={1\over 2\,{\tt m_1}}\,{\rm L}^2+O(\LL)
\end{equation}
\begin{equation}\label{varas}
{\rm var\,} \widehat{K}_{\rho}= {{\tt m}_2\over 3\,{\tt m}_1^3}\,{\rm L}^3 +O({\rm L}^{2})\,.
\end{equation}

\end{theorem}

\vskip0.5cm
\noindent{\bf Justification of the decay condition} 
Our plan is to decompose $g_1=h_1+h_2$ so that $h_1$ will absorb the logarithmic part of $g_1$ at $\infty$
and will have a manageable Mellin transform, while $M h_2$ will satisfy the decay condition by  virtue of the 
following classical result.

\vskip0.5cm

\begin{lemma}\label{class} {\rm \cite[Section 1.29]{titch}} Let $\phi$ be analytic in a sector $-\alpha<{\rm Arg}\,\rho<\beta$ with $0<\alpha,\beta\leq\pi$, and
suppose that in this sector for some $a<b$ and each $\delta>0$
$$\phi(\rho)=O(|\rho|^{-a-\delta}) \quad{\rm as~~}|\rho|\downarrow 0\,,~\phi(\rho)=O(|\rho|^{-b+\delta}) 
\quad{\rm as~~}|\rho|\uparrow\infty\,.$$
Then $M\phi$ is analytic in the strip $a<\Re\,s<b$ and satisfies
$$M\phi(s)=O\left(e^{-(\beta-\epsilon)\Im s}\right)\,\quad {\rm as~} \Im s\to\infty\,,\quad
M\phi(s)=O\left(e^{(\alpha-\epsilon)\Im s}\right)\,\quad {\rm as~} \Im s\to -\infty\,,\qquad$$
for each $\epsilon>0$ uniformly in every strip strictly inside $a<\Re\, s<b$.
\end{lemma}

\par By definition
$$g_1(\rho)=\int_0^1 f^{(1)}(\rho(1-x))(1-e^{-\rho x})\tOm({\rm d}x)$$
which is an entire function, where 
$$f^{(1)}(\rho)=e^{-\rho}\sum_{n=1}^\infty {\rho^n\over n!}\,{\mathbb E\,}K_n$$
and the series in $\rho$ has all coefficients positive. The same applies to the series involved in 
$$1-e^{-\rho}=e^{-\rho}\sum_{j=1}^\infty {\rho^j\over j!}.$$
Therefore $|g_1(\rho)|\leq g_1(|\rho|)$ and by (\ref{2g})
$$|g_1(\rho)|=O(\LL^3)\,~~~|\rho|\to\infty\,.$$
Writing (\ref{2g}) as
$$g_1(\rho)=\sum_{j=0}^3 d_j\LL^j+O(\rho^{-\epsilon})\,,\quad \rho\to\infty\,,~\rho\in {\mathbb R}$$
and selecting
$$h_1(\rho):=d_0e^{-1/\rho}+\log (\rho+1)\sum_{j=1}^3 d_j \log^{j-1}\rho$$
we obtain a function $h_2:=g_1-h_1$ which is analytic in the open halfplane $\Re\,\rho>0$,
has the same expansion in powers of $\LL$ as $g_1$ for $\rho\to\infty$ and
satisfies
Lemma \ref{class} with $a=1$, $b=\epsilon$ and $\alpha=\beta=\pi/2$. Hence by the lemma
$$|M h_2(s)|=O\left(e^{-(\pi/2-\delta)|s|}\right)\,,\quad {\rm for~~}-1/2<\Re s<\epsilon\,,~|s|\to\infty\,$$ 
for each $\delta$. 
The meromorphic continuation 
of $Mh_1$ from $-1<\Re\, s<0$ to the halfplane $\Re\,s> -1$ follows 
from the standard Mellin transform formulas
$$M e^{-1/\,\cdot\,}(s)=\Gamma(-s), ~~M\log(1+\cdot)(s)=-\Gamma(s)\Gamma(-s)\,,~~
M(\phi(\,\cdot\,)\log(\,\cdot\,))(s)={{\rm d}\over {\rm d}s}M\phi(s)$$
which by application of (\ref{stirling}) also imply 
$$|M h_1(s)|=O\left(e^{-(\pi/2-\delta)|s|}\right)\,,\quad {\rm for~~}-1/2<\Re s<\epsilon\,,~|s|\to\infty\,.$$ 
 
\par It follows that in a strip containing the imaginary axis $|Mg_1(s)|=O\left(e^{-(\pi/2-\delta)|s|}\right)$, 
and the same estimate holds for 
$|Mf^{(2)}(s)|=|Mg_1(s)/\Phi(-s)|$ due to Lemma \ref{fourier}.
Thus  $|Mf^{(2)}|$ decays exponentially fast and  Lemma \ref{singan} is applicable.

\vskip0.5cm
\noindent
{\bf Number of parts prior to the exponential split} Let $\tau$ be an independent exponential time, with rate $\lambda$.
By (\ref{f-split-rec}), the Mellin transform in $\rho$ satisfies 
$$M({\mathbb E}\widehat{K}_{\rho}(\tau))(s)={M \widehat{\Phi}(s)\over \lambda+\Phi(-s)}.$$
The term $\lambda$ in the denominator kills the zero of $\Phi$ at $s=0$, therefore arguing as in the analysis of
 (\ref{Mf1}) 
we see that this expression has only {\it double} pole at $s=0$.
It follows that 
$${\mathbb E}\widehat{K}_{\rho}(\tau)= d\,{\rm L} + O(1)$$
for each $\lambda>0$, with $d$ depending on $\lambda$.
Repeatedly appealing to (\ref{f-split-rec}) and Lemmas \ref{fourier}, \ref{le2} and \ref{class}  we obtain
$f^{(m)(\rho)}=O(\LL^m)$, which trivially implies the following rough estimate

\begin{lemma}\label{splitmom}
For $\tau$ an independent exponential time,
as $\rho \te \infty$ with ${\rm L } = \log \rho$,
$${\mathbb E}(\widehat{K}_\rho(\tau))^m=O({\rm L}^m)\,,\qquad m=1,2,\ldots$$
\end{lemma}
\vskip0.5cm
\noindent
{\bf Depoissonisation} Modifying the argument in \cite[Section 6.2]{gpy03regvar} we obtain for 
$\rho=n\to \infty$
$$\widehat{K}_n\sim K_n\,, \,~~\widehat{K}_n(t)\sim K_n(t)$$
for each $t>0$, almost surely, with all factorial moments, and with at least two central moments.
The same sandwich-type proof works because the moments grow logarithmically
(thus vary regularly as $n\to\infty$).

\section{Central limit theorem for $K_n$}

We switch to fixed-$n$ framework, in order to apply 
a delicate recent result due to Neininger and R{\"u}schendorf, which we reproduce below
with a minor adaptation and  notational changes.
Suppose a sequence $(Y_n)$ of random variables satisfies
\begin{equation}\label{distr-eqn}
Y_n\stackrel{d}{=}b_n+Y_{I_n}'\,,\qquad n\geq 1
\end{equation}
where $(Y_n)\stackrel{d}{=}(Y_n')$, the pair
$(b_n, I_n)$ is independent of $(Y_n')$, each $I_n$ assumes values in $\{0,\ldots, n\}$, 
and ${\mathbb P}(I_n=n)<1$ for
$n\geq 1$. Let $\mu_n={\mathbb E} \,Y_n$, $\sigma_n^2={\rm var\,}Y_n$, and $\parallel\cdot\parallel_3$ denote the
${\cal L}^3$-norm of a random variable.
\vskip0.5cm

\begin{theorem}\label{nr} {\rm \cite[Theorem 2.1]{neinruesch}}
Assume that each $\parallel Y_n\parallel_3<\infty$ and that $(Y_n)$ satisfies the recursion \rem{distr-eqn}.
Suppose that for some constants $C>0,\,\alpha>0$ the following three
conditions all hold:

\begin{enumerate}

\item[{\rm (i)}] 
$$
\lim\sup_{n\to\infty} {\mathbb E}\,\log\left({I_n\vee 1
\over n}\right)<0\,,\qquad \sup_{n>1} \bigg\vert\bigg\vert \log\left({I_n\vee 1\over n } \right) \bigg\vert\bigg\vert_3<\infty$$

\item[{\rm (ii)}] for some $\lambda\in[0,\,2\alpha[\, $ and some $\kappa$
$$\parallel b_n-\mu_n +\mu_{I_n} \parallel_3=O(\log^{\kappa}n)\,,\qquad \sigma_n^2= C\,\log^{2\alpha}n+O(\log^{\lambda}n) $$

\item[{\rm (iii)}] 
$$\alpha>{1\over 3}+\,\max\left(\kappa\,,\,{\lambda\over 2}\right) .$$
\end{enumerate}
Then the law of 
$${Y_n-\mu_n\over \sqrt{C} \log^{\alpha} n}$$
converges weakly to the standard normal distribution.

\end{theorem}
\vskip0.5cm
\noindent
Condition (i) says that the split $I_n/n$ must be bounded away from $0$ and $1$,
and the further two conditions  require logarithmic growth of moments. See \cite{neinruesch} for proof and 
estimates of the distance between the probability laws.
\vskip0.5cm

\par To apply the result to $K_n$
we split the range of subordinator at an independent exponential time $\tau$ with mean $1$.
Let $b_n=K_n(\tau)$ be the number of parts of ${\cal C}_n$ produced by the multiplicative subordinator 
in the period between $0$ to $\tau$.
Observe that 
(\ref{distr-eqn}) holds, with $I_n$ being the number of uniform sample points larger than $\widetilde{S}_{\tau}$.
The asymptotics of moments in Theorem \ref{momenty} and Lemma \ref{splitmom}
apply to $K_n$ and $b_n$ literally without change, 
by the virtue of depoissonisation,

\par 
Let us check now the conditions of the theorem.
Since $1\leq K_n\leq n$, all absolute moments of $K_n$ are finite. Furthermore,
the variable $I_n/n$ converges strongly and with all moments to $1-\widetilde{S}_{\tau}$,
thus $-\log (I_n/n)$ approaches the stopped value $S_{\tau}$ of the {\it additive} subordinator,
which is positive 
and
has ${\mathbb E} S_{\tau}^3<\infty $ as a consequence of ${\tt m}_3<\infty$,
by an easy application of 
the L{\'e}vy-Khintchine
formula. Whence (i) is satisfied. 
In evaluating the asymptotics of moments we can switch to the poissonised composition with $\rho=n$, ${\rm L}=\log n$.
Then from Theorem \ref{momenty} we see that 
$\parallel\mu_n-\mu_{I_n}\parallel_3$ is of the order of $\rm L$.
By Lemma \ref{splitmom} also $\parallel b_n\parallel_3=O({\rm L})$, because 
${\mathbb E}\,b_n^3= O({\rm L}^3)$.
By the above and Theorem \ref{momenty} the parameters involved in (ii) are $\alpha=3/2$, $\lambda=2$ and $\kappa=1$,
so condition (iii) is satisfied.
Thus all  conditions of Theorem \ref{nr} are fulfilled and we 
deduce the following conclusion:

\vskip0.5cm
\begin{theorem}
\label{clt} If the L{\'e}vy measure has a continuous density on $]0,\infty[$ and satisfies 
{\rm (L)} and {\rm (R)}  then the distribution of the random variable
$${K_n- \,{\rm L}^2/(2{\tt m_1})\over \sqrt{{{\tt m}_2/ (3\,{\tt m}_1^3})}\,\,{\rm L}^{3/2}}\,,\qquad {\rm L}=\log n$$
converges weakly to the standard normal distribution, as $n\to\infty$.
\end{theorem} 

\vskip0.5cm
\noindent

\section{Small parts of the composition}

We have shown in \cite{gpy03regvar} that
in the regular variation case all $K_{n,r}$, $r=1,2,\ldots$ are of the same order of growth as $K_n$ and, suitably normalised,
converge almost surely to multiples of the same random variable.
In the compound Poisson case these variables are bounded as $n$ grows \cite{gnedin03bs}, and
if the increments of $(S_t)$ are exponentially distributed (this corresponds to the Ewens partition structure)
the $K_{n,r}$'s converge to independent Poisson random variables \cite{abt}. 
 Our next goal is proving a joint central limit theorem for the small-part counts.

\vskip0.5cm
\noindent
\paragraph{Marginal central limit theorems} The line of argument  repeats that for $K_n$.
Let $f^{(m)}_r(\rho)$ be the $m$th factorial moment of $\widehat{K}_{\rho,r}$, e.g.
$f^{(1)}_1(\rho)$ is the mean number of singletons. Recall that the recursion (\ref{f-rec})
holds for $f^{(m)}=f^{(m)}_r$
with $\pi(\rho)= {e^{-\rho }\rho^r/ r!}$.
Applying the Mellin transform to the recursion we obtain
\begin{equation}\label{Mfr}
M f_r^{(1)}(s)= {\Gamma(r+s)\over r!} {\Phi(-s:-s)\over \Phi(-s)}\,,\qquad -1 < \Re \,s < 0.
\end{equation}
This agrees with (\ref{Mf1}) and $K_\rho=\Sigma_r \, K_{\rho,r}$ in view of the identity
$$-\Gamma(s)=\sum_{r=1}^\infty {\Gamma(r+s)\over \Gamma(r+1)}\,,~~~~~-1 < \Re \,s < 0.$$ 
The right-hand side of (\ref{Mfr}) can be extended meromorphically through the imaginary axis, with a {\it double} pole at $0$, where
we have the Laurent expansion

\begin{equation}\label{LMr}
M f_r^{(1)}(s)= {1\over {\tt m}_1 r}\,{s^{-2}}-
d_1\,{s^{-1}}+O(1)\,
\end{equation}
where
$$d_1=
{2c{\tt m}_1-2\gamma\,{\tt m}_1+
{\tt m}_2-2{\tt m}_1\,\psi^{(0)}(r)
\over 2{\tt m}_1^2\, r}$$
and $\psi^{(0)}$ is the digamma function. The decay at complex infinity is justified as before, thus 
by Lemma \ref{singan}
\begin{equation}\label{f1r}
f^{(1)}_r(\rho) = {{\rm L}\over {\tt m}_1 r}+d_1 + O(\rho^{-\epsilon}).
\end{equation}

\par  Now we need to translate (\ref{LMr}) into asymptotics of the right-hand side of (\ref{f-rec}) with
$m=2$. We note that evaluation of integrals 
with a factor of $e^{-\rho x}x^r$ essentially 
amounts to Tauberian-type asymptotics, since for $r>1$ the factor $x^r$ compensates the singularity of $\tOm$ at $x=0$, and 
$e^{-\rho x}$ is negligible outside each fixed vicinity of $0$. In particular, we have

\vskip0.5cm
\begin{lemma} If the condition {\rm \,\,(L)} holds then for $\rho\to\infty$
$$ \int_0^1 {e^{-\rho x}(\rho x)^r\over r!}\,\tOm({\rm d}x)={1\over r}+O(\rho^{-\epsilon})\,.$$
\end{lemma}
{\it Proof.} Integrating by parts and replacing $\vec{\nu}$ by the right-hand side of 
(L), and pushing the upper integration limit to $\infty$
the claim is reduced to the standard integral
$$\int_0^{\infty} (-\log x+d) 
\left({e^{-\rho x}\rho^r x^{r-1}\over (r-1)!}-{e^{-\rho x}\rho^{r+1} x^r\over r!}\right){\rm d}x={1\over r}\,$$
which does not depend on $d$, because 
$$\int_0^{\infty} 
\left({e^{-\rho x}\rho^r x^{r-1}\over (r-1)!}-{e^{-\rho x}\rho^{r+1} x^r\over r!}\right){\rm d}x=0\,
$$
as is easily checked. \endpf
\vskip0.5cm
\par Using the lemma and (\ref{f1r}) we compute the right-hand side of (\ref{f-rec}) as
$$g(\rho) = {2{\rm L}\over {\tt m}_1 r^2}+{2d_1\over r} + O(1),$$
hence by Lemma \ref{singan}
$$Mg(s) = {2\over {\tt m}_1 r^2}\,{s^{-2}}-{2d_1\over r}\,{s^{-1}} + O(1)\,.$$
We compute then
$$Mf^{(2)}_r(s)={Mg(s)\over \Phi(-s)} = -{2\over {\tt m}_1^2r^2}\,s^{-3}+\left({{\tt m}_2\over r^2{\tt m_1}^3}+
{2d_1\over r{\tt m}_1}\right)\,s^{-2} + O(s^{-1})$$
which yields again by Lemma \ref{singan} 
$$f^{(2)} _r (\rho) = {1\over r^2{\tt m}_1^2}{\rm L}^2+ \left({{\tt m}_2\over r^2{\tt m_1}^3}+
{2d_1\over r{\tt m}_1}\right)
{\rm L} + O(1)$$
(the decay condition is again checked by application of Lemma \ref{fourier}).
This together with (\ref{f1r}) implies 
$${\rm var}\,K_{r,\rho} = \left({{\tt m}_2\over r^2{\tt m}_1^3}+{1\over r{\tt m}_1}\right)
{\rm L} + O(1) \,.
$$

\par Decomposing as in (\ref{fpoint}) at rate $1$ exponential time, the Mellin transform of the expectation of
$b_\rho=\widehat{K}_{\rho,r}(\tau)$ has $\Phi(-s)+1$ in the denominator, hence
all moments of
$b_\rho$ remain bounded as $\rho\to\infty$. 
Passing to the fixed-$n$ 
version and applying Theorem \ref{nr} 
with $\alpha=1/2,\lambda=\kappa=0$ we obtain
\vskip0.5cm

\begin{theorem} 
\label{CLT-univ}
Under our assumptions on the L{\'e}vy measure,
for large $n$ the distribution of each $K_{n,r}$ 
is approximately normal, with moments
$${\mathbb E}K_{n,r} = {1\over {\tt m}_1 r}\log n\, + O(1),\qquad 
{\rm var}\,K_{n,r} =
\left({{\tt m}_2\over r^2{\tt m}_1^3}+{1\over r{\tt m}_1}\right)
\log n + O(1)\,.$$
The same is true for $\widehat{K}_{\rho,r}$ as $\rho\to\infty$.
\end{theorem}

\vskip0.5cm
\noindent
{\bf Joint central limit theorem} Our strategy to prove a joint central limit theorem for the small parts counts
is to consider 
more general finite patterns and 
functionals of composition such as $\Sigma_r \,a_r \widehat{K}_{\rho,r}$. By linearity, each variable of this kind decomposes 
as in (\ref{fpoint}), hence the method we applied to $\widehat{K}_{\rho,r}$ can be used to justify the normal limit.
As an instance of such a functional consider
$$\widehat{K}_{\rho, [r]}:=\widehat{K}_{\rho,1}+\cdots+\widehat{K}_{\rho,r}$$
and define $f^{(m)}_{[r]}(\rho)$ to be the $m$th factorial moment of
$\widehat{K}_{\rho, [r]}$. The corresponding pattern is $\{1,\ldots,r\}$, thus the moments
satisfy (\ref{f-rec}) with $\pi(\rho)=\sum_{j=1}^r e^{-\rho} {\rho^j/j!}$.
Obviously from (\ref{f1r})

$$f^{(1)}_{[r]} = {h_r\over{\tt m}_1}{\rm L}+ O(1)$$
where
$h_r:=\sum_{j=1}^r 1/j$ are the  harmonic numbers.
Letting Mellin's machine roll to produce $f_{[r]}^{(1)}\to g\to Mg\to M f_{[r]}^{(2)}\to f_{[r]}^{(2)}$
the variance is computed as 
$${\rm var}
\,K_{\rho, [r]} = \left({{\tt m}_2 h_r^2\over {\tt m}_1^3}+{h_r\over {\tt m}_1}\right)
{\rm L} + O(1).
$$
\par Similar computation for the pattern $E=\{i,j\}$ and Theorem \ref{CLT-univ} yield the covariance 
$$
{\rm cov}(\widehat{K}_{\rho,i},\widehat{K}_{\rho,j}) = \left({{\tt m}_2\over{\tt m}_1^3}{1\over i\,j}+1(i=j){1\over j\,{\tt m}_1}
\right)\,{\rm L} + O(1)\,.
$$
\vskip0.5cm

\begin{theorem} Under our assumptions on $\nu$, as $n\to\infty$, 
the infinite random sequence
$$\big((K_{n,r}-{\mathbb E}K_{n,r})\log^{-1/2} n\,,\,\,\,\, r=1,2,\ldots,n\big)$$
 converges in law to a multivariate gaussian 
sequence with the covariance matrix
\begin{equation}\label{cov}
\left({{\tt m}_2\over{\tt m}_1^3}{1\over i\,j}+1(i=j){1\over j\,{\tt m}_1}
\right)_{i,j=1}^{\infty}\,.
\end{equation}
\end{theorem}
{\it Proof.} To justify the joint gaussian law it suffices to establish convergence 
to a gaussian limit for each finite linear combination
$$Q_n=\sum_r a_r K_{n,r}\,.$$
From the above facts about $K_{n,r}$ and from (\ref{cov}) it follows that 
both the expectation and the variance of $Q_n$ grow like $\rm L$.
Also, due to an obvious additivity, $Q_n$ 
satisfies a distributional equation of the type (\ref{distr-eqn}) provided we split
the range of subordinator by the value at an independent exponential time $\tau$.
Then $b_n$ is equal to the contribution to $Q_n$ by the jumps of subordinator before $\tau$.
The conditions of Theorem \ref{nr} are checked exactly as for $K_n$ or $K_{n,r}$, 
thus by this theorem $Q_n$ is approximately gaussian.
\endpf

\section{The gamma case and further examples}

In the case of  the gamma subordinators we have

$$\quad c=-\log\theta\,,\qquad{\tt m}_j={(j-1)!\over \theta^j}\,,$$
and the singular expansion 
$$\Phi(-s:-s)=-s^{-1}+(c-\gamma)-d_2 s+O(s^2)$$
with constant
$$d_2=\int_0^1 {\log x\left. x(1-x)^{\theta-1} +\log(1-x)\over -x\log(1-x)\right.}\,{\rm d}x$$
which does not seem to simplify.
The Mellin transform of ${\mathbb E}\widehat{K}_\rho$ is given by the formula
$$M f^{(1)}(s)={-\Gamma(s)\over \log(1-s/\theta)}\int_0^1{x^{-s}(1-x)^{\theta-1}\over -\log(1-x)}\,{\rm d}x$$ 
which defines a function which is meromorphic for ${\Re\,}s>-1$ 
with a sole singularity at 
$s=0$. The Laurent expansion of $M f^{(1)}$ at $s=0$ involves
$${\tt m}_j={(j-1)!\over \theta^j}\,,\quad c=-\log\theta\,,\quad d_1={\gamma^2\over 2}+{\pi^2\over 12}\,.$$
Leaving only principal terms, the 
asymptotics of moments is
$${\mathbb E}\widehat{K}_{\rho} = {\theta\over 2}\,{\rm L}^2 + O(\LL)\,,\qquad
{\rm var\,}\widehat{K}_{\rho}= {\theta\over 3}\,{\rm L}^3 + O(\LL^2)
$$
as was promised in Theorem \ref{gammaK}.
\vskip0.5cm

\noindent
{\bf Further examples} Another instance of a gamma-type subordinator
was introduced in \cite{gnedin03tsf} to describe a composition 
resembling the ordered Ewens sampling formula. This subordinator 
has (multiplicative) L{\'e}vy measure 
$$\tOm({\rm d}x)=(1-x)^{\theta-1}x^{-1}\,{\rm d}x\,,\qquad x\in \,]0,1]$$
with parameter $\theta>0$. 
The Laplace exponent given by the formula
$$
\Phi(s)=\sum_{j=1}^\infty \left({1\over j+\theta-1}-{1\over j+\theta-1 +s}\right)\,
$$
is a function interpolating the generalised harmonic numbers
\begin{equation}\label{harmtheta}
\Phi(n)=\sum_{j=1}^n {1\over j+\theta-1}\,,
\end{equation}
which can be seen as a combinatorial analogue of the logarithmical Laplace exponent for the gamma subordinator.
\par The basic characteristics of subordinator are readily expressed in terms of the polygamma function
$$\psi^{(k)}(\theta)={{\rm d}^{k+1}\, \log \Gamma (\theta)\over {\rm d}\,\theta^{k+1}}.$$
Thus we have
$$\vec{\nu}(x)=-\log x-\psi^{(0)}(\theta)+O(x^{-1})$$
as is shown by expanding
$$\vec{\nu}(x)=\int_x^1(1-z)^{\theta-1}z^{-1}\,{\rm d}z=-\log x+c(\theta)-\gamma+c_1(\theta)x+\ldots\,,$$
substituting $\theta=1$ to see that $c(1)=\gamma=-\psi^{(0)}(1)$, then differentiating in $\theta$ and sending $x\to 0$
to obtain
\begin{equation}\label{c'}
c'(\theta)=\int_0^1 (1-z)^{\theta-1}z^{-1}\log(1-z)\,{\rm d}z=-\psi^{(1)}(\theta)\,,
\end{equation}
whence $c(\theta)=-\psi^{(0)}(\theta)$.
It follows as in Appendix or directly from (\ref{harmtheta}) that the expansion at infinity is 
$$\Phi(\rho)=\LL -\psi^{(0)}(\theta)+O(\rho^{-1})\,.$$
The moments are computed by further differentiating (\ref{c'}) as
$$
{\tt m}_j=\int_0^1 (1-x)^{\theta-1}x^{-1}|\log(1-x)|^j\,{\rm d}x=
(-1)^{j+1} \psi^{(j)}(\theta).$$

\par See \cite{berg} for computation of some densities related to this family of subordinators,
and \cite{py03h} for further examples of L{\'e}vy measures with logarithmic singularity.

\vskip0.5cm

\noindent
{\bf Oscillatory asymptotics for a discrete measure} The constant term $c$ assumed in (L) cancels
in the asymptotics of the moments. However, this assumption is essential by our approach.
Let us assume  
a bounded oscillating term in place of $c$, and examine 
how the asymptotics could be affected.
For ease of computation we consider
the atomic measure 
$$\tOm({\rm d}x)=\sum_{j=1}^{\infty} \delta_{1/e^{j}}({\rm d}x)\,.$$
For this measure $\Phi(n)$ is equal to the expected maximum in a sample of $n$ geometric random
variables, which was analysed in 
\cite[Example 12]{flaj-harm}.

\par Denoting $\lfloor\,\cdot\,\rfloor$, respectively $\{\,\cdot\,\}$, the integer and the fractional parts of a positive number,
we have 
$$\vec{\nu}(x)=\lfloor -\log x\rfloor=-\log x-\{-\log x\}\,,~~~~x\in\,]0,1]$$
thus there is a logarithmic singularity but
the expansion (L) does not hold, because the second term oscillates between $0$ and $1$.
Condition (R) is satisfied since $\vec{\nu}(x)=0$ for $x>1/e$ and all moments are finite,
$${\tt m}_k=\sum_{j=1}^\infty |\log(1-e^{-j})|^k\,,$$
e.g. ${\tt m}_1=0.6843\,,~{\tt m}_2=0.2345$ (truncated at four decimals).
Furthermore, 
$$\Phit (\rho)=\sum_{j=1}^{\infty}(1-e^{-\rho/e^k})\,,\quad M \Phit(s)={-\Gamma(s)\over 1-e^s}\,,\quad
\Phi(s)=\sum_{j=1}^{\infty} (1-(1-e^{-j})^s).$$

\par It is seen that $M\Phit$ has a {\it double} pole at $s=0$ with Laurent expansion
$$M\Phit(s) ={1\over s^2}-{\gamma+1/2\over s}+{1\over 2}(1+6\gamma+6\gamma^2+\pi^2)+O(s)$$
and infinitely many {\it simple} poles on the imaginary axis at 
$$s_k=2\pi{\tt i}k\,,\quad k\in {\mathbb Z}\setminus\{0\}.$$ 
The conclusion of Lemma \ref{singan} still holds, which can be  justified by applying the inverse Mellin 
transform formula
and integrating over increasing rectangular contours, as $k\to\infty$, with sides $\Im s =2\pi k+1/2,~\Im s =-2\pi k-1/2,~
\Re s=1/2,~\Re s=-1/2$. Thus

\begin{equation}\label{tiny}
\Phit(\rho)=\LL+\gamma+1/2+\phi(\LL)+O(\rho^{-1+\epsilon})\,,\qquad \rho\to\infty
\end{equation}
where the contribution of the imaginary poles amounts to the term $\phi(\LL)$ given by the formula
$$\phi(u)=-\sum_{k\in{\mathbb Z}\setminus \{0\}}\Gamma(2\pi {\tt i}k)\,e^{2\pi{\tt i}k u}\,,$$
which is a periodic function with period one and a tiny amplitude. The fluctuations of the $O(1)$ term in (\ref{tiny}) 
are not asymptotically  negligible, though they are very small due to the fast decay of the Fourier coefficients,
e.g. 
$$\Gamma(s_1)=\overline{\Gamma(s_{-1})}=(0.126+0.501\,{\tt i})\,10^{-4}.$$
(truncated at seven decimals).

\par It follows that
$Mf^{(1)}$ has a triple pole at $0$ and  simple imaginary poles $s_k=2\pi {\tt i}k$, hence 
$$f^{(1)}={1\over 2{\tt m}_1}\LL^2+\left({\gamma+1/2\over {\tt m}_1}+{{\tt m}_2\over 2{\tt m}_1^2}\right)\LL+
\phi_1(\LL)+O(\rho^{-1+\epsilon})$$
where $\phi_1$ is another periodic function.
Thus oscillation prevails only in the third term in the expansion of $f^{(1)}$. However,
computing $g_1$, as in Lemma \ref{le2}, 
we should include $\LL^2\Phit(\rho)$, thus arriving at
$$g_1(\rho)={1\over {\tt m}_1}\LL^3+ \phi_2(\LL)\LL^2+\ldots$$
with oscillating $\phi_2$. 
As a consequence we will have 
$$f^{(2)}= {1\over 4{\tt m}_1^2}\LL^4+\phi_3(\LL)\LL^3+\ldots,$$
and this suggests that the principal $O(\LL^3)$-term in the asymptotic expansion of the variance 
${\rm var\,} \widehat{K}_\rho$ will oscillate,
though we could not establish this fact rigorously.

\par In a similar situation of sampling from the geometric distribution, the expectation of the 
number of different values in a sample
 also involves an 
oscillating second term, and the same is true for the factorial moment of order 2 \cite{prodinger}. 
However, this has no impact on the principal asymptotics of the variance:
the oscillating terms cancel and the variance 
converges to a constant, as had been shown long ago 
by Karlin, see \cite[Example 6, p. 385]{karlin67urn}.

\section{Appendix}

\begin{lemma} The four conditions in {\rm (L)} are equivalent. \end{lemma}
{\it Proof.} 
The equivalence of expansions of $\nu$ and $\tOm$ follows by the change of variables
$x=1-e^{-y}$ and $y=-\log(1-x)$. 
For example, assuming the expansion of $\nu$ we obtain that of $\tOm$ by
using 
$$\tOm[x,1]=\nu [-\log(1-x),\infty]$$
and substituting
$$-\log(-\log(1-x))= -\log (x(1+x/2+\ldots))=-\log x+O(x).$$
\par The expansions of integrals are equivalent to the expansions of measures by writing
$$\Phi(\rho)/\rho=\int_0^{\infty} e^{-\rho y}\nu[y,\infty]\,{\rm d}y\,,\qquad
\Phit(\rho)/\rho=\int_0^{1} e^{-\rho x}\tOm[x,1]\,{\rm d}x\,$$
and using 
the classical integral
$$\rho \int_0^{\infty} e^{-\rho y} \log (1/y)\,{\rm d}y={\rm L}+\gamma\,$$
together with standard properties of the Laplace transform. 
\endpf

\begin{lemma} The function 
$$\Phi(s)=\int_0^1 (1-(1-x)^s) \,\tOm({\rm d}x)$$
has no zeros for ${\rm Re}\, s>0$. And if $\Phi(s)=0$ for purely imaginary $s$, then 
$\tOm$ is atomic, with support $\{1-a^k, k=1,2,\ldots\}$ for some $0<a<1$.

\end{lemma}
{\it Proof.} For ${\Re}\, s>0$ and $x\in \,]0,1[$
we have $|(1-x)^s|<1$. Therefore ${\Re}\, (1-(1-x)^s)>0$ and the integral cannot be zero.
For real $r\neq 0$ the equality
$\Phi({\rm i}\, r)=0$ is only possible when $1=(1-x)^{{\rm i}\,r}$ holds $\tOm$-almost everywhere,
but such $x$ is of the form $x=1-\exp(-2\pi k/|r|)$ for some $k=1,2,\ldots$.
\endpf

\def\cprime{$'$} \def\polhk#1{\setbox0=\hbox{#1}{\ooalign{\hidewidth
\lower1.5ex\hbox{`}\hidewidth\crcr\unhbox0}}} \def\cprime{$'$}
\def\cprime{$'$} \def\cprime{$'$}
\def\polhk#1{\setbox0=\hbox{#1}{\ooalign{\hidewidth
\lower1.5ex\hbox{`}\hidewidth\crcr\unhbox0}}} \def\cprime{$'$}
\def\cprime{$'$} \def\polhk#1{\setbox0=\hbox{#1}{\ooalign{\hidewidth
\lower1.5ex\hbox{`}\hidewidth\crcr\unhbox0}}} \def\cprime{$'$}
\def\cprime{$'$} \def\cydot{\leavevmode\raise.4ex\hbox{.}} \def\cprime{$'$}
\def\cprime{$'$} \def\cprime{$'$} \def\cprime{$'$}


\begin{thebibliography}{10}


\bibitem{Armstrong} J.A.~Armstrong and N.~Bleistein,
\newblock Asymptotic expansions of integrals with oscillatory kernels and logarithmic singularities,
\newblock {\it SIAM J. Math. Anal.}, 11: 300--307, 1980. 


\bibitem{abt} R.~Arratia, A.G.~Barbour and S. Tavar{\'e}
\newblock {\em Logarithmic combinatorial structures: A probabilistic approach}
\newblock 2002 (forthcoming book)


\bibitem{berg} C.~Berg and A.J.~Duran
\newblock Some transformations of Hausdorff moment sequences and harmonic numbers, 
\newblock preprint, 2004.


\bibitem{Fedoryuk} M.V.~Fedoryuk, 
\newblock {\it Asymptotics: Integrals and series.}
\newblock Nauka, Moscow, 1987.





\bibitem{flaj-harm}
P.~Flajolet, X.~Gourdon and P.~Dumas,
\newblock Mellin transforms and asymptotics: harmonic sums.
\newblock {\em Theoret. Comput. Sci.} 144: 3--58, 1995.


\bibitem{gnedin97} A.V.~ Gnedin, 
The representation of composition structures,
\newblock {\it Ann. Probab.} 25: 1437-1450, 1997.


\bibitem{gnedin03bs}
A.V.~Gnedin.
\newblock The Bernoulli sieve,
\newblock {\em Bernoulli} 10: 79-96, 2004.


\bibitem{gnedin03tsf}
A.V.~Gnedin.
\newblock {Three sampling formulas},
\newblock {\em Combinatorics, Probability and Computing} 13: 185-193, 2004.


\bibitem{gnedinp03}
A.V.~Gnedin and J.~Pitman.
\newblock {Regenerative composition structures}.
\newblock {\it Ann. Prob.} (to appear)


\bibitem{gpy03regvar}
A.V.~Gnedin, J.~Pitman and M.~Yor.
\newblock {Asymptotic laws for compositions derived from transformed subordinators}.
\newblock (available at arXiv:math.PR/0403438)


\bibitem{karlin67urn}
S.~Karlin.
\newblock Central limit theorems for certain infinite urn schemes.
\newblock {\em J. Math. Mech.}, 17:373--401, 1967.


\bibitem{neinruesch}
R.~Neininger and L.~R{\"u}schendorf, 
\newblock On the contraction method with degenerate limit equation.
\newblock {\it Ann. Prob.} (to appear)


\bibitem{csp}
J.~Pitman.
\newblock {Combinatorial stochastic processes}.
\newblock Technical Report 621, Dept. Statistics, U.C. Berkeley, 2002.
\newblock Lecture notes for St. Flour course, July 2002. (available via
www.stat.berkeley.edu)


\bibitem{prodinger} H.~Prodinger.
\newblock Compositions and patricia tries: no fluctuations in the variance! preprint 2003.



\bibitem{py03h}
J.~Pitman and M.~Yor.
\newblock {Infinitely divisible laws associated with hyperbolic functions}.
\newblock {\em Canad. J. Math.}, 55(581):292--330, 2003.

\bibitem{titch} E.C.~Titchmarsh
\newblock {\it Introduction to the theory of Fourier integrals,}
\newblock Oxford Univ. Press 1937.
 

\end{thebibliography}
\end {document}